\documentclass[11pt, equation]{article}
\usepackage{amssymb}
\usepackage{amsmath}
\usepackage{theorem}
\evensidemargin\oddsidemargin
\makeatletter
\@addtoreset{equation}{section}

\makeatother
\newtheorem{theo}{Theorem}

\newtheorem{lem}{Lemma}

{\theorembodyfont{\rmfamily}  }
\theoremheaderfont{\scshape}
\newcommand{\PP}{\Bbb{P}}

\newcommand{\R}{\Bbb{R}}

\newcommand{\E}{\Bbb{E}}

\newcommand{\be}{\begin{equation}}
\newcommand{\ee}{\end{equation}}
\newcommand{\fc}{\mathcal{F}}
\numberwithin{theo}{section}
\numberwithin{equation}{section}
\numberwithin{propo}{section}
\numberwithin{lem}{section}
\numberwithin{defi}{section}

{\theorembodyfont{\rmfamily}  

\begin{document}
\date{}
\title{A general Darling-Erd\H os theorem in Euclidean space}
\author{ Gauthier Dierickx\thanks{Research supported by a PhD fellowship of the Research Foundation Flanders-FWO } \;\;and  Uwe Einmahl\\
Department of Mathematics, Vrije Universiteit Brussel\\
Pleinlaan 2, B-1050 Brussels, Belgium}
\maketitle

\begin{abstract}
We provide an improved version of the Darling-Erd\H os theorem for sums of i.i.d. random variables with mean zero and finite variance. We extend this result to multidimensional random vectors. Our proof is based on a new strong invariance principle in this setting which has other applications as well such as an integral test refinement of the multidimensional  Hartman-Wintner LIL. We also identify a borderline situation where one has weak convergence to a shifted version of the standard limiting distribution in the classical Darling-Erd\H os theorem.
\end{abstract}


\noindent \textit{AMS 2010 Subject Classifications:} 60F15, 60F17.

\noindent \textit{Keywords:} Darling-Erd\H os theorem, extreme value distribution, Hartman-Wintner LIL, integral test, strong invariance principle, multidimensional version, double truncation.

\section{Introduction} Let $X, X_1, X_2, \ldots$ be i.i.d. random variables and set $S_n=\sum_{j=1}^n X_j, n \ge1.$  Further set $Lt =\log_e(t \vee e), LLt = L(Lt)$ and $LLLt = L(LLt), t \ge 0.$  In 1956 Darling and Erd\H os proved that under the assumption $\E|X|^3 < \infty,$ $\E X^2=1$ and $\E X=0,$ the following convergence in distribution result holds,
\be \label{1956}
a_n \max_{1\le k \le n} |S_k|/\sqrt{k} - b_n \stackrel{d}{\to}  \tilde{Y},
\ee
where $a_n = \sqrt{2LLn}, b_n = 2 LLn + LLLn/2 - \log(\pi)/2$ and $\tilde{Y}$ is a random variable which has an extreme value distribution with distribution function $y \mapsto \exp(-\exp(-y)).$\\
The above third moment assumption was later relaxed in \cite{ood} and \cite{sho}  to $\E |X|^{2 + \delta} < \infty$ for some $\delta >0,$ but the question remained open whether a finite second moment would be already sufficient.\\
This was finally answered in \cite{E-5}, where it is shown that (\ref{1956}) holds if and only if 
$$\E X^2 I\{|X| \ge t\}= o ((LLt)^{-1}) \mbox{ as }t \to \infty.$$
Moreover, it is shown in \cite{E-5} that the above result holds more generally under the assumption of a finite second moment if one replaces the normalizers $\sqrt{k}$ by $\sqrt{B_k}$, where $B_n=\sum_{j=1}^n \sigma_j^2$ and
$$\sigma_n^2:= \E X^2 I\{|X| \le \sqrt{n}/(LLn)^p\}$$
for some $p \ge 2.$ So we have under the classical assumption $\E X^2 < \infty$ and $\E X=0,$
 $$a_n \max_{1\le k \le n} |S_k|/\sqrt{B_k} - b_n \stackrel{d}{\to}  \tilde{Y}.$$
 For some further related work on the classical Darling-Erd\H os theorem the reader is referred to \cite{be},\cite{csw},\cite{EM} and the references in these articles.\\
 The Darling-Erd\H os theorem  is also  related to finding an integral test refining the Hartman-Wintner LIL, a problem which was already addressed by Feller in 1946.  Here one can  relatively easily prove that the classical Kolmogorov-Erd\H os-Petrowski integral test for Brownian motion holds for sums of i.i.d. mean zero random variables  if one has $\E |X|^{2 + \delta} < \infty$ for some $\delta >0$. In this case one has for any non-decreasing function $\phi: ]0,\infty[ \to ]0,\infty[,$
 $$\PP\{ |S_n| \le \sqrt{n}\phi(n)\mbox{ eventually}\} = 1 \mbox { or } 0,$$
 according as
 $\sum_{n=1}^{\infty} n^{-1}\phi(n)\exp(-\phi^2(n)/2)$ is finite or infinite. \\
 Feller proved that this result remains valid under the second moment assumption if one  replaces $\sqrt{n}$ by $\sqrt{B_n}$ defined as above with $p \ge 4$.  Similarly as for the Darling-Erd\H os theorem this implies that the Kolmogorov-Erd\H os-Petrowski integral test holds in its original form if  
 $$\E X^2 I\{|X| \ge t\}= O ((LLt)^{-1}) \mbox{ as }t \to \infty.$$
  The proof in \cite{Fe} was based on a skillful double truncation argument which only worked for symmetric distributions. Finally in \cite{E-5} an extension of this argument to the general non-symmetric case was found so that we now know that most results in \cite{Fe} are correct. (See also \cite{bai} for more historical background.)\\
 There is still one question in the paper \cite{Fe} which has not yet been addressed, namely whether it 
 is possible to make the theorem ``slightly more elegant'' by replacing the sequence $\sqrt{B_n}$ by $\sqrt{n}\sigma_n.$ Feller writes that he ``was unable to devise a proof simple enough to be justified by the slight improvement of the theorem'' (see p. 632 in \cite{Fe}). We believe that we have found  a simple enough proof  of Feller's claim. (See Step 3 in the proof of Theorem \ref{asip1}.) \\
This  leads to the following improved version of the Darling-Erd\H os theorem under the finite second moment assumption:
 \be
 a_n \max_{1\le k \le n} |S_k|/\sqrt{k} \sigma_k - b_n \stackrel{d}{\to}  \tilde{Y}.
 \ee
 At the same time we can show that there is a much wider choice for the truncation level in the definition of $\sigma_n^2$. For instance, it  is possible to define $\sigma_n^2$ as $\E X^2I\{|X| \le \sqrt{n}\}.$\\
This improved version of the Darling-Erd\H os theorem will actually follow from a general result for $d$-dimensional random vectors which will be given in the following section.

\section{Statement of Main Results} 
We now consider i.i.d. $d$-dimensional random vectors $X, X_1, X_2, \ldots$  such that $\E|X|^2 < \infty$ and $\E X=0$, where we denote the Euclidean norm by $|\cdot|$. The corresponding matrix norm will be denoted by $\|\cdot\|$, that is, we set
$$\|A\|:=\sup_{|x| \le 1} |Ax|$$ 
for any $(d,d)$-matrix $A$. It is well known that $\|A\|=$ the largest eigenvalue of $A$ if   $A$ is symmetric and non-negative definite.\\Let  again $S_n:=\sum_{j=1}^n X_j, n \ge 1.$   Horv\'ath \cite{Hor} obtained in 1994 the following multidimensional version of the Darling-Erd\H os theorem assuming that $\E |X|^{2+\delta} < \infty$ for some $\delta > 0$ and that  Cov($X$) (=  the covariance matrix of $X$) is equal to the $d$-dimensional identity matrix $I$,
\be \label{Hor}
a_n \max_{1\le k \le n} |S_k|/\sqrt{k} - b_{d,n} \stackrel{d}{\to} \tilde{Y},
\ee
where $\tilde{Y}$ has the same distribution as in dimension 1, 
$$b_{d,n}:= 2LLn + d LLLn/2 - \log(\Gamma(d/2)),$$
and $\Gamma(t), t > 0$ is the Gamma function. Recall that $\Gamma(1/2)=\sqrt{\pi}$ so that this extends the 1-dimensional Darling-Erd\H os theorem. \\[.1cm]
We are ready to formulate our general result. We consider non-decreasing sequences $c_n$ of positive real numbers satisfying for large $n$,
\be \label{c_n}
\exp(-(\log n)^{\epsilon_n}) \le c_n/\sqrt{n} \le \exp((\log n)^{\epsilon_n}),
\ee
where $\epsilon_n \to 0.$\\
Further, let for each $n$, $\Gamma_n$ be the symmetric non-negative definite matrix such that 
\be  \label{gamma_n}
\Gamma_n^2 = [\E X^{(i)}X^{(j)} I\{|X| \le c_n\}]_{1 \le i, j \le d}, n \ge 1.
\ee
If the covariance matrix of $X=(X^{(1)},\ldots, X^{(d)})$ is positive definite, the matrices $\Gamma_n$ will be invertible for large enough $n$. Replacing $c_n$ by $c_{n \vee n_0}$ for a suitable $n_0 \ge 1$ if necessary, we can assume w.l.o.g.  that {\bf all} matrices $\Gamma_n$ are invertible.
\begin{theo} \label{DE-general}
Let $X, X_n, n \ge 1$ be i.i.d. mean zero random vectors in $\R^d$ with $\E |X|^2 < \infty$ and $\mathrm{Cov}(X) = I.$ Then we have for any sequence $\{c_n\}$ satisfying condition (\ref{c_n}),
\be \label{DEtheo}
a_n \max_{1 \le k \le n} |\Gamma_k^{-1} S_k|/\sqrt{k} - b_{d,n} \stackrel{d}{\to} \tilde{Y},
\ee
where $\tilde{Y}: \Omega \to \R$ is a random variable such that 
$$\PP\{\tilde{Y} \le t\}=\exp(-\exp(-t)), t \in \R.$$ Under the additional assumption
\be \label{log-log}
\E|X|^2 I\{|X| \ge t\}=o((LLt)^{-1}) \mbox{ as } t \to \infty,
\ee
we also have 
\be \label{DEtheo1}
a_n \max_{1 \le k \le n} |S_k|/\sqrt{k} - b_{d,n} \stackrel{d}{\to} \tilde{Y}.
\ee
\end{theo}
It is easy to see that condition (\ref{log-log}) is satisfied if $\E |X|^2LL|X|< \infty.$ This latter condition, however, is more restrictive than (\ref{log-log}).\\ It is natural to ask whether this condition is also necessary as in the 1-dimensional case. (See Theorem 2 in \cite{E-5}.) This question becomes much more involved in the multidimensional case and  we get a slightly weaker result, namely that the following condition
\be \label{fconv}
\E|X|^2I\{|X| \ge t\}= O((LLt)^{-1}) \mbox{ as } t \to \infty
\ee
is necessary for (\ref{DEtheo1}). \\
To prove this result we show that if  condition (\ref{fconv}) is not satisfied,  then \\$a_n \max_{1 \le k \le n} |S_k|/\sqrt{k} - b_{d,n}$  cannot converge in distribution  to any variable of the form  $\tilde{Y}+c,$ where $c \in \R.$\\ If one allows this larger class of limiting distributions,  condition (\ref{fconv}) is optimal. There are examples where $\E|X|^2I\{|X| \ge t\} = O((LLt)^{-1})$ and\\ $a_n \max_{1 \le k \le n} |S_k|/\sqrt{k} - b_{d,n}$ converges in distribution to   $\tilde{Y}+c$ for some $c \ne 0.$ (See Theorem \ref{shift} below.)
\begin{theo} \label{converse}
 Let $X, X_n, n \ge 1$ be i.i.d. mean zero random vectors in $\R^d$ with $\E |X|^2 < \infty$ and $\mathrm{Cov}(X) = I$ and suppose that there exists a $c \in \R$ such that
 \be \label{EXTRAc}
 a_n \max_{1 \le k \le n} |S_k|/\sqrt{k} - b_{d,n} \stackrel{d}{\to} \tilde{Y} + c,
 \ee
 where $\tilde{Y}$ is as in Theorem \ref{DE-general}. Then condition (\ref{fconv}) holds. \end{theo}
Our basic tool for proving the above results is a new strong invariance principle for sums of i.i.d. random vectors which is valid under a finite second moment assumption. If one has an approximation with an almost sure error term of order $o(\sqrt{n/LLn})$, one can obtain the Darling-Erd\H os theorem directly from the normally distributed case. The problem is that it is impossible to get such an approximation under the sole assumption of a finite second moment.  In \cite{E-5} it was shown that one needs the ``good'' approximation of order $o(\sqrt{n/LLn})$ only if the sums $|S_n|$ are large and it was shown that  in dimension 1  one can obtain approximations which are particularly efficient for the random subsequence where the sums are large.  Using recent results on $d$-dimensional strong approximations (see \cite{sakh} and \cite {E-2}), we are now able to obtain an analogue of the approximation in \cite{E-5} in the $d$-dimensional setting (see Lemma \ref{Sterke_Approx} and  relation (\ref{step3f1}) below). As an additional new feature we also show that an approximation by $\Gamma_n \sum_{j=1}^n Z_j$ is possible where $Z_j, j\ge 1$ are i.i.d. $\mathcal{N}(0,I)$-distributed random vectors with $\mathcal{N}(0,\Sigma)$ denoting the $d$-dimensional normal distribution with mean zero and covariance matrix $\Sigma.$ This type of  approximation leads to the improved versions of the Darling-Erd\H os theorem  and Feller's integral test as indicated in Section 1. 
\begin{theo} \label{asip1} Let $X, X_n, n \ge 1$ be i.i.d. mean zero random vectors in $\R^d$ with $\E |X|^2 < \infty$ and $\mathrm{Cov}(X) = \Gamma^2,$ where $\Gamma$ is a symmetric non-negative definite (d,d)-matrix. 
Let  $c_n$ be a non-decreasing sequence of positive real-numbers satisfying condition (\ref{c_n}) for large $n$ and
let $\Gamma_n$ be defined as in (\ref{gamma_n}).
 If the underlying p-space $(\Omega, \fc, \PP)$ is rich enough one can construct independent $\mathcal{N}(0,I)$-distributed random vectors $Z_n, n \ge 1$ such that we have for the partial sums $T_n:=\sum_{j=1}^n Z_j, n \ge 1,$
\begin{itemize}
\item [(a)] $S_n - \Gamma \;T_n = o(\sqrt{n LLn})$ as $n \to \infty$ with prob. 1, 
\item [(b)] $\PP\{ |S_n - \Gamma_n T_n| \ge 2\sqrt{n}/ LLn, |S_n| \ge \frac{4}{3}\|\Gamma\|\sqrt{n LLn} \mbox{ infinitely often}\}=0,$ and
\item [(c)] $\PP\{ |S_n - \Gamma_n T_n| \ge 2\sqrt{n}/LLn, |\Gamma T_n| \ge \frac{4}{3}\|\Gamma\|\sqrt{n LLn} \mbox{ infinitely often}\}=0.$
\end{itemize}
 \end{theo}
 Combining our strong invariance principle with the Kolmogorov-Erd\H os-Petrowski integral test for $d$-dimensional Brownian motion, one obtains by the same arguments as in Section 5 of \cite{E-5} the following result,\begin{theo} 
 \label{test} Let $X, X_n, n \ge 1$ be i.i.d. mean zero random vectors in $\R^d$ with $\E |X|^2 < \infty$ and $\mathrm{Cov}(X) = \Gamma^2,$ where $\Gamma$ is a symmetric positive definite (d,d)-matrix. 
 Let  $c_n$ be a non-decreasing sequence of positive real-numbers satisfying condition (\ref{c_n}) for large $n$ and let $\Gamma_n$ be defined as in (\ref{gamma_n}).
 Then we have for any non-decreasing function $\phi: ]0,\infty[ \to ]0,\infty[,$
 $$\PP\{|\Gamma_n^{-1}S_n| \le \sqrt{n}\phi(n) \mbox{ eventually}\} = 1 \mbox{ or }=0$$
 according as 
 $$I_d(\phi):= \sum_{n=1}^{\infty}n^{-1} \phi^d(n)\exp(-\phi^2(n)/2) < \infty \mbox{ or} =\infty.$$
 \end{theo}
 Note that we can assume w.l.o.g. that all matrices $\Gamma_n$ are invertible since they converge to $\Gamma$ which is invertible. \\
Let $\lambda_n (\Lambda_n)$ be the smallest (largest) eigenvalue of $\Gamma_n, n \ge 1.$ Assuming that Cov($X$) $=I$, we can infer from Theorem \ref{test},
 \be
 I_d(\phi) < \infty \Rightarrow \PP\{|S_n| \le \Lambda_n \sqrt{n}\phi(n) \mbox{ eventually}\} = 1
\ee
 and
 \be
 I_d(\phi) = \infty \Rightarrow \PP\{|S_n| >\lambda_n \sqrt{n} \phi(n) \mbox{ infinitely often}\} = 1,
\ee
which is the $d$-dimensional version of the result conjectured by Feller in \cite{Fe}.\\[.2cm]
The proof of our strong invariance principle (= Theorem \ref{asip1}) will be given in Sect. 3. In the two subsequent sections 4 and 5 we will show how Theorems \ref{DE-general} and \ref{converse} follow from the strong invariance principle. In Sect. 6 we return to the real-valued case and show that if $\E X^2I\{|X| \ge t\} \sim c(LLt)^{-1}$ that then (\ref{1956})  still remains valid if we replace $\tilde{Y}$ by $\tilde{Y}-c$. Finally, we answer a question which was posed in \cite{kls}.
 
 \section{Proof of the strong invariance principle}
 Our proof is divided into three steps.\\[.1cm]
 STEP 1. We recall a double truncation argument which goes back to Feller \cite{Fe} for symmetric random variables. This was later extended to non-symmetric random variables in \cite{E-5}  and finally to random elements in Hilbert space  in \cite{E-3}. To formulate the relevant result we need some extra notation. We set
 \begin{align*}
	&X_n^{'} := X_n I \{ | X_n | \leq \sqrt{n}/ (LLn)^5 \} , && \overline{X}'_n := X'_n - \mathbb{E}X'_n;\\
	&X''_n := X_n I \{  \sqrt{n}/ (LLn)^5 < | X_n | \leq \sqrt{n LLn} \} , && \overline{X}''_n := X''_n - \mathbb{E}X''_n;\\
	&X'''_n := X_n I \{    \sqrt{n LLn} < | X_n |  \} , && \overline{X}'''_n := X'''_n - \mathbb{E}X'''_n;
	\end{align*}
and we denote the corresponding sums by $S'_n, \overline{S}'_n,  S''_n, \overline{S}''_n, S'''_n, \overline{S}'''_n.$\\
Then we have (see \cite{E-3}, Lemma 11 and Lemma 12) 
\be \label{step1f1}
 S_n - \overline{S}'_n = o(\sqrt{nLLn}) \mbox{ a.s.}
 \ee
and
\be \label{step1f2}
\PP\{ |S_n - \overline{S}'_n| \ge \sqrt{n}/ (LLn), |S_n| \ge \|\Gamma\|\sqrt{n LLn} \mbox{ i.o.}\} =0.
 \ee
STEP 2. Let $\Sigma_n$ be the sequence of symmetric non-negative definite matrices such that $\Sigma_n^2$ is the covariance matrix of $X'_n$ for $n \ge 1.$ Furthermore,  let $A(t)$ be the  symmetric non-negative definite matrices satisfying  $$A(t)^2 = [\E X^{(j)}X^{(k)}I\{|X| \le t\}]_{1 \le j, k \le d}, t \ge 0.$$  It is easy to see that $A(t)^2, t \ge 0$  is monotone, that is,  $A(t)^2 - A(s)^2$ is non-negative definite if $0 \le s \le t.$\\ This implies that  $A(t), t \ge 0$ is monotone as well (see Theorem V.1.9 in \cite{Bha}). Consequently, $A(c_n), n \ge 1$ is a monotone sequence of symmetric non-negative definite matrices whenever $c_n$ is non-decreasing. Moreover, $A(c_n)$ converges to $\Gamma$ if $c_n \to \infty.$ \\[.2cm]  We  have the following strong approximation result, where we set 
$$\tilde{\Gamma}_n:=A(\sqrt{n}/(LLn)^5), n \ge 1.$$
\begin{lem}\label{Sterke_Approx}
 If the underlying p-space is rich enough,
	one can construct independent random vectors $Z_i \sim \mathcal{N}(0, I)$ such that	
	\be \label{eq2-1}
	\overline{S}'_n -  \sum_{i=1}^n \tilde{\Gamma}_i Z_i = o(\sqrt{n}/LLn ) \text{ a.s. }\ee
	and
	\be \label{eq-2}
	S_n -  \sum_{i=1}^n \Gamma Z_i = o(\sqrt{nLLn}) \mbox{ a.s.}
	\ee
\end{lem}
{\bf Proof}
(i) We first show that one can construct independent $\mathcal{N}(0,I)$-distributed random vectors such that
$$ \overline{S}'_n -  \sum_{i=1}^n \Sigma_i Z_i = o(\sqrt{n}/LLn ) \text{ a.s. }$$
By Corollary $3.2$ from \cite{E-2} and the fact that $\mathbb{E} |  \overline{X}'_n  |^3 \leq 8  \mathbb{E} |  X'_n |^3,$
  it is enough to show 
	$$\sum_{n =1}^{\infty} \mathbb{E} \left|  X'_n  \right|^3 / \left(  \frac{ \sqrt{ n }  }{LLn }  \right)^3 < \infty.$$
	Using the simple inequality, 
	$$ \mathbb{E} |  X'_n |^3 \le  \mathbb{E} |X|^{2 + \delta} I\{|X| \le \sqrt{n}\}  \sqrt{n}^{1-\delta}(LLn)^{-5(1-\delta)}, 0 < \delta <1,$$ we find (setting $\delta=2/5$) that the above series is 
	$$\le \sum_{n=1}^{\infty} \mathbb{E} |X|^{12/5} I\{|X| \le \sqrt{n}\} /\sqrt{n}^{12/5}$$
	Using a standard argument (see, for instance, the proof of part (a) of Lemma 3.3 in \cite{E-2}), one can show that this last series is finite whenever $\E|X|^2 < \infty.$ \\[.1cm]
(ii) To complete the proof of (\ref{eq2-1}) it is now sufficient to show that
$$\sum_{i=1}^n (\Sigma_i - \tilde{\Gamma}_i)Z_i = o(\sqrt{n}/LLn) \mbox{ a.s.}$$
 By a standard argument this follows if 
 \be \label{series}
 \sum_{n=1}^{\infty} \frac{  \mathbb{E} \left| \left( \Sigma_n - \tilde{\Gamma}_n  \right) Z_n  \right|^2  
 }{ n/(LLn)^2 } < \infty.
 \ee
 Since $Z_n \sim \mathcal{N}(0, I)$, we have $$\mathbb{E} \left| \left( \Sigma_n - \tilde{\Gamma}_n  \right) Z_n  \right|^2  
 \leq 
  d  \left\| \Sigma_n - \tilde{\Gamma}_n  \right\|^2 \le d\|\Sigma_n^2 -\tilde{ \Gamma}_n^2\|,   $$
  where we have used Theorem X.1.1 in \cite{Bha} for the last inequality. 
 From the definition of $\Sigma_n$ and $\tilde{\Gamma}_n$ it is obvious that
 $$\langle x, (\tilde{\Gamma}_n ^2-\Sigma_n^2 )x \rangle  = \left(\mathbb{E}  \langle X, x\rangle  I \{ | X| \leq  \sqrt{n}/(LLn)^5 \}  \right)^2,  x \in \mathbb{R}^d.$$
 The last expression equals $ \left(\mathbb{E} \langle X, x\rangle  I \{ | X| >  \sqrt{n}/(LLn)^5 \}  \right)^2$ since $ \mathbb{E} \langle X, x\rangle =0 .$\\
 Hence $\left\| \Sigma_n^2 - \tilde{\Gamma}_n^2  \right\| \leq \left(\mathbb{E}  | X|   I \{ | X| >  \sqrt{n}/(LLn)^5 \} \right)^2 \le (\E|X|^2)^2 n^{-1} (LLn)^{10} $. \\
 It is easy now to see that the series in (\ref{series}) is finite.\\[.1cm]
 (iii) Finally note that 
 $$S_n - \sum_{i=1}^n \Gamma Z_i = (S_n -  \overline{S}'_n) +  (\overline{S}'_n -\sum_{i=1}^n \tilde{\Gamma}_i Z_i) + \sum_{i=1}^n (\tilde{\Gamma}_i -\Gamma)Z_i,$$
 where the first two terms are of almost sure order $o(\sqrt{nLLn})$ by (\ref{step1f1})  and  (\ref{eq2-1}), respectively. Since $\tilde{\Gamma}_n \to \Gamma$ as $ n \to \infty$, we also have that
 $$\sum_{i=1}^n (\tilde{\Gamma}_i -\Gamma)Z_i=o(\sqrt{nLLn}) \mbox{ a.s.},$$
 and we can conclude that indeed $S_n - \sum_{i=1}^n \Gamma Z_i = o(\sqrt{nLLn})$ a.s.
 Lemma \ref{Sterke_Approx} has been proven. 
 $\Box$\\[.2cm]
 STEP 3. Combining Lemma \ref{Sterke_Approx} with relations (\ref{step1f1}) and  (\ref{step1f2}) we find that 
 \be \label{step3f1}
 \PP\left\{|S_n - \sum_{j=1}^n \tilde{\Gamma}_j Z_j | \ge \ 3\sqrt{n}/(2 LLn), |\Gamma T_n| \ge \frac{5}{4}\|\Gamma\|\sqrt{n LLn} \mbox{ i.o.}\right\} =0.
 \ee
We  next show that
\be \label{step3f2}
\PP\left\{ |\sum_{j=1}^n (\Gamma_n - \tilde{\Gamma}_j)Z_j | \ge \sqrt{n}/(2 LLn), 	|T_n| \ge \frac{5}{4}\sqrt{n LLn} \mbox{ i.o.}\right\} =0,
\ee
where $$\Gamma_n:=A(c_n), n \ge 1$$ and $c_n$ is an arbitrary non-decreasing sequence of positive real numbers satisfying condition (\ref{c_n}) for large $n.$\\[.1cm]
 Using  that $\{|\Gamma T_n| \ge \|\Gamma\| x\} \subset \{|T_n| \ge x\}, x >0$, we get from   (\ref{step3f1}) and  (\ref{step3f2}):
\be \label{step4f1}
\PP\{|S_n - \Gamma_n T_n | \ge 2\sqrt{n}/ LLn, |\Gamma T_n| \ge \frac{5}{4}\|\Gamma\|\sqrt{n LLn} \mbox{ i.o.}\} =0.
\ee
Further recall that $S_n -\Gamma T_n = o(\sqrt{nLLn})$ a.s. (see Lemma \ref{Sterke_Approx}).  Consequently, we can infer from 
(\ref{step4f1}) that
\be \label{step4f2}
\PP\{|S_n - \Gamma_n T_n | \ge 2\sqrt{n}/ LLn, |S_n| \ge\frac{4}{3} \|\Gamma\|\sqrt{n LLn} \mbox{ i.o.}\} =0.
\ee
We see that  the proof of Theorem \ref{asip1} is complete once we have established (\ref{step3f2}). Toward this end
we need the following inequality which is valid for normally distributed random vectors $Y: \Omega \to \R^d$ with mean zero and covariance matrix $\Sigma$:
\be \label{step3f3}
\PP\{|Y| \ge x\} \le \exp(-x^2/(8\sigma^2)), x \ge 2\E |Y|^2,	
\ee
where $\sigma^2$ is the largest eigenvalue of $\Sigma.$ (See Lemma 4 in \cite{E-3}.)\\
From (\ref{step3f3}) we trivially get  that
\be \label{step3f4}
\PP\{|Y| \ge x\} \le 2\exp(-x^2/(8\E |Y|^2)), x \ge 0.
\ee
Though this last inequality is clearly suboptimal, it will  nevertheless be more than sufficient for the proof of (\ref{step3f2}).\\[.2cm]
{\bf Proof of (\ref{step3f2}).} To simplify notation we set $d_n:=\sqrt{n}/(LLn)^5, n \ge 1,$ $n_k:=2^k$ and $\ell_k:=[2^{k-1}/(Lk)^5], k \ge 0.$ By the Borel-Cantelli lemma it is enough to show that 
$$\sum_{k=1}^{\infty}\PP\left(\bigcup_{n=n_{k-1}}^{n_k}\left\{ |\sum_{j=1}^n (\Gamma_n - \tilde{\Gamma}_j)Z_j | \ge \sqrt{n}/(2 LLn), |T_n| \ge \frac{5}{4} \sqrt{n LLn} \right\} \right) < \infty.$$
Set $\tilde{\mathbb{N}}:=  \tilde{\mathbb{N}}_1 \cap \tilde{\mathbb{N}}_2$, where
$$ \tilde{\mathbb{N}}_1:=\{k: \|\Gamma_{n_k} - \tilde{\Gamma}_{\ell_k}\| \le \|\Gamma\| (Lk)^{-5/2}\}$$
and
$$\tilde{\mathbb{N}}_2:= \{k: \|\Gamma_{n_k} - \Gamma_{n_{k-1}}\| \le (Lk)^{-2}\}.$$ Then it is easy to see that the above series is finite if 
\be \label{step3f5}
\sum_{k \in \tilde{\mathbb{N}}} \PP\left\{ \max_{n_{k-1} \le n \le n_k}  |\sum_{j=1}^n (\Gamma_n -\tilde{ \Gamma}_j)Z_j | \ge 2^{(k-1)/2}/(2 Lk)\right\} < \infty
\ee
and
\be \label{step3f6}
\sum_{k \not \in \tilde{\mathbb{N}}} \PP\left\{ \max_{n_{k-1} \le n \le n_k}|T_n| \ge  2^{(k-1)/2} (Lk)^{1/2}\right\} < \infty.
\ee
To bound the series in (\ref{step3f5}), we  first note that employing the L\'evy inequality for sums of  independent symmetric random vectors, one obtains
\begin{eqnarray*}
&&\PP\left\{ \max_{n_{k-1} \le n \le n_k}  |\sum_{j=1}^n (\Gamma_n - \tilde{\Gamma}_j)Z_j | \ge 2^{(k-1)/2}/(2 Lk)\right\} \\&\le& \PP\left\{ \max_{n_{k-1} \le n \le n_k}  |\sum_{j=1}^n (\Gamma_{n_{k}} - \tilde{\Gamma}_j)Z_j | \ge 2^{(k-1)/2}/(4 Lk)\right\}\\&&\hspace{.3cm}  +\; \PP\left\{ \max_{n_{k-1} \le n \le n_k}  |(\Gamma_{n_k} - \Gamma_n)\sum_{j=1}^n Z_j | \ge 2^{(k-1)/2}/( 4Lk)\right\} \\
&\le& 2 \PP\left\{   |\sum_{j=1}^{n_k} (\Gamma_{n_{k}} - \tilde{\Gamma}_j)Z_j | \ge 2^{(k-1)/2}/(4Lk)\right\}
\\&&\hspace{.3cm}
+\;2\PP\left\{\|\Gamma_{n_k} - \Gamma_{n_{k-1}}\|\;|\sum_{j=1}^{n_k} Z_j| \ge 2^{(k-1)/2}/(4 Lk)\right\}=:2 p_{k,1} + 2 p_{k,2},
\end{eqnarray*}
where we have also used the fact that $$\|\Gamma_{n_k}-\Gamma_n\| \le \|\Gamma_{n_k}-\Gamma_{n_{k-1}}\| \mbox{ if } n_{k-1} \le n \le n_k.$$ This follows easily from the monotonicity of the sequence $\Gamma_n.$\\[.1cm]
To bound $p_{k,1}$, we first note that by Theorem X.1.1 in \cite{Bha} for $1 \le j \le \ell_k,$
\begin{eqnarray}
\|\Gamma_{n_k} - \tilde{\Gamma}_j\|^2 &\le& \|\Gamma^2_{n_k} - \tilde{\Gamma}^2_j\| = \sup_{|x| \le 1} \E\langle x, X \rangle^2 I\{d_j \wedge c_{n_k} < |X| \le d_j \vee c_{n_k}\} \nonumber \\
&\le&\sup_{|x| \le 1} \E\langle x, X \rangle^2=\|\Gamma\|^2. \label{step3f6a}
 \end{eqnarray}
Apply (\ref{step3f4}) with $Y= \sum_{j=1}^{n_k} (\Gamma_{n_{k}} - \tilde{\Gamma}_j)Z_j.$ Then clearly $\E Y = 0$ and, moreover, by independence of the random vectors $Z_j$ we have for $k \in \tilde{\mathbb{N}},$
\begin{eqnarray*}
\E |Y|^2 &=&\sum_{j=1}^{n_k} \E |(\Gamma_{n_{k}} - \tilde{\Gamma}_j)Z_j|^2 \le d \sum_{j=1}^{n_k} \|\Gamma_{n_k} - \tilde{\Gamma}_j\|^2\\&\le& d\|\Gamma\|^2(\ell_k + (n_k-\ell_k)(Lk)^{-5})\le 2d\|\Gamma\|^2 n_k(Lk)^{-5}.
\end{eqnarray*}
We  conclude that
$$p_{k,1} \le 2\exp(-(Lk)^3/(512d\|\Gamma\|^2)), k \in \tilde{\mathbb{N}}_1.$$
Similarly, we obtain
$$p_{k,2} \le \PP\left\{|\sum_{j=1}^{n_k} Z_j| \ge 2^{(k-1)/2} Lk/4\right\} \le 
2\exp(-(Lk)^2/(256d)), k \in \tilde{\mathbb{N}}_2.$$
It is now clear that the series in (\ref{step3f5}) is finite.\\[.1cm]
To show that the series in (\ref{step3f6}) is finite, we note that by (\ref{step3f4}) and the L\'evy inequality,
$$ \PP\left\{ \max_{n_{k-1} \le n \le n_k}|T_n| \ge  2^{(k-1)/2} (Lk)^{1/2}\right\} 
\le 2 \PP\left\{|T_{n_k}| \ge  2^{(k-1)/2} (Lk)^{1/2}\right\} \le 4 k^{-\eta},$$
where $\eta=(16d)^{-1}$ and it is enough to check that
$$\sum_{k \not \in \tilde{\mathbb{N}}}k^{-\eta} < \infty.$$
To verify that this series is finite, observe that by the argument used in (\ref{step3f6a}) we have,
\begin{eqnarray*}
\|\Gamma_{n_k}-\Gamma_{n_{k-1}}\|^2 &\le & \sup_{|x| \le 1}\E\langle x, X \rangle^2 I\{ c_{n_{k-1}} < |X| \le \ c_{n_k}\}\\ &\le& \E|X|^2 I\{ c_{n_{k-1}} < |X| \le \ c_{n_k}\},
\end{eqnarray*}
which implies 
$$\sum_{k=1}^{\infty} \|\Gamma_{n_k}-\Gamma_{n_{k-1}}\|^2 \le \E |X|^2 < \infty.$$
We conclude that
$$
\sum_{k \not \in \tilde{\mathbb{N}}_2}(Lk)^{-2} < \infty.$$
So the proof of (\ref{step3f2}) is complete if we show that
\be \label{extra}
 \sum_{k \not \in \tilde{\mathbb{N}}_1}k^{-\eta} < \infty.
 \ee
We need another lemma.
\begin{lem}
Consider two sequences $c_{k,i}, k \ge 1$ of positive real numbers satisfying for large enough $k,$
\be \label{ass_lem32}
2^{k/2}\exp(-k^{\delta}) \le c_{k,i} \le 2^{k/2}\exp(k^{\delta}) ,i=1,2,
\ee
where $0<\delta <1.$ Set $\Gamma_{k,i} := A(c_{k,i}), i=1,2, k \ge 1.$
Then we have,
$$\sum_{k=1}^{\infty}k^{-\delta} \|\Gamma_{k,1} - \Gamma_{k,2}\|^2 < \infty.$$
\end{lem}
{\bf Proof.} Using the same argument as in (\ref{step3f6a}), we have for large $k,$
$$
\|\Gamma_{k,1} - \Gamma_{k,2}\|^2 \le \E |X|^2 I\{2^{k/2}\exp(-k^{\delta}) <|X| \le 2^{k/2}\exp(k^{\delta})\} \le \sum_{j=[k- 3k^{\delta}]}^{ [k+ 3k^{\delta}]} \beta_j,
$$
where $\beta_j :=\E|X|^2 I\{2^{j-1} < |X|^2 \le 2^j\}, j \ge 1$.\\ We can conclude that for some $k_0 \ge1$ and a suitable $j_0 \ge 0,$
\be \label{step3f12}
\sum_{k=k_0}^{\infty}k^{-\delta} \|\Gamma_{k,1} - \Gamma_{k,2}\|^2  \le \sum_{j=j_0}^{\infty}\beta_j \sum_{k=m_1(j)}^{m_2(j)} k^{-\delta},
\ee
where $$m_1(j)=\min\{k \ge k_0: [k+3k^{\delta}] \ge j\}$$ and  $$m_2(j)=\max\{k \ge k_0: [k-3k^{\delta}] \le j\}.$$  It is easy to see that $m_1(j) \ge j-3j^{\delta} \ge j/2$ and $m_2(j) \le j + 4j^{\delta}$ for large $j$.\\
Consequently, we have for large $j,$
\be  \label{step3f11}
\sum_{k=m_1(j)}^{m_2(j)} k^{-\delta} \le 2^{\delta} (m_2(j)-m_1(j)+ 1)j^{-\delta} \le 2^{3+\delta} < \infty.
\ee
We obviously have $\sum_{j=1}^{\infty} \beta_j < \infty$ (as $\E |X|^2 < \infty$). Combining relations (\ref{step3f12}) and  (\ref{step3f11}) we obtain the assertion of the lemma. $\Box$\\[.1cm]
We apply the above lemma with $c_{k,1} = c_{n_k}, c_{k,2} = d_{\ell_k}, k \ge 1.$ From condition (\ref{c_n}) we readily obtain that for large $k$, $$2^{k/2}\exp(-k^{\epsilon'_k}) \le c_{k,1} \le 2^{k/2}\exp(k^{\epsilon'_k}),$$ where $\epsilon'_k:=\epsilon_{n_k} \to 0$ so that condition (\ref{ass_lem32}) is satisfied for any $\delta > 0$. This is also the case for the sequence $c_{k,2}$. So we can choose $\delta=\eta/2$ and it follows that
$$\infty >\sum_{k \not \in \tilde{\mathbb{N}}_1}k^{-\eta/2}\|\Gamma_{n_k} -\tilde{\Gamma}_{\ell_k}\|^2 \ge  \|\Gamma\|^2\sum_{k \not \in \tilde{\mathbb{N}}_1} k^{-\eta/2}(Lk)^{-5},$$
which shows that (\ref{extra}) holds.
 \section{Proof of Theorem \ref{DE-general}}
We first prove (\ref{DEtheo}).
Set $k_n= [\exp((Ln)^{\alpha})]$, where $0 < \alpha <1$. Then it follows from the $d$-dimensional version of the Hartman-Wintner LIL that  for any given $\epsilon >0,$ with prob. 1, 
$$|\Gamma_k^{-1} S_k|/\sqrt{k} \le \lambda_k^{-1}\sqrt{2LLk}(1+\epsilon), k \ge k_0(\omega, \epsilon),$$
where $\lambda_k$ is the smallest eigenvalue of $\Gamma_k$. As $\lambda_k \nearrow 1,$
we can conclude that for large enough $n,$
$$\max_{1 \le k <k_n} |\Gamma_k^{-1} S_k|/\sqrt{k} \le \sqrt{2\alpha LLn}(1+2\epsilon),$$
which is $\le \sqrt{2LLn}$ if we choose $\epsilon$ small enough.
It follows that
\be \label{proof1}
a_n \max_{1 \le k <k_n} |\Gamma_k^{-1} S_k|/\sqrt{k} -b_{d,n} \to -\infty \mbox{ a.s.}
\ee
So (\ref{DEtheo}) holds if and only if  
$$a_n \max_{k \in K_n} |\Gamma_k^{-1} S_k|/\sqrt{k} -b_{d,n} \stackrel{d}{\to} \tilde{Y},$$
where $K_n:=\{k_n+1,\ldots, n\}.$\\
We split $K_n$  into two random subsets:
$$K_{n,1}(\cdot):=\{k \in K_n: |S_k - \Gamma_k T_k| \le 2\sqrt{k}/(LLk)\},
K_{n,2}(\cdot):= K_n \setminus K_{n,1}(\cdot).$$
In view of Theorem \ref{asip1}(b) (where we set $\Gamma=I$) there are with prob. 1 only finitely many $k$'s such that
 $$|S_k - \Gamma_k T_k| > 2\sqrt{k}/LLk \mbox{ and }|\Gamma_k^{-1} S_k| \ge 4 \sqrt{kLLk}/(3\lambda_k),$$ where $\lambda_k$ is  again the smallest eigenvalue of $\Gamma_k$. \\As $\lambda_k \nearrow 1,$ we can conclude that with prob. 1 there are only finitely many $k$'s such that 
 $$|S_k - \Gamma_k T_k| > 2\sqrt{k}/LLk \mbox{ and }|\Gamma_k^{-1} S_k| \ge  \sqrt{2kLLk},$$ and it follows that
$$a_n \max_{k \in K_{n,2}(\cdot)} |\Gamma_k^{-1} S_k|/\sqrt{k} -b_{d,n} \to -\infty \mbox{ a.s.}$$
We see that (\ref{DEtheo}) is equivalent to  
$$a_n \max_{k \in K_{n,1}(\cdot)} |\Gamma_k^{-1} S_k|/\sqrt{k} -b_{d,n} \stackrel{d}{\to} \tilde{Y}.$$
From the definition of the sets $K_{n,1}(\cdot)$ we easily get that
$$ a_n \max_{k \in K_{n,1}(\cdot)} |\Gamma_k^{-1} S_k|/\sqrt{k} - a_n \max_{k \in K_{n,1}(\cdot)} | T_k|/\sqrt{k} \to 0 \mbox{ a.s.}$$
By Slutsky's lemma  (\ref{DEtheo}) holds if and only if  
$$a_n \max_{k \in K_{n,1}(\cdot)} |T_k|/\sqrt{k} -b_{d,n} \stackrel{d}{\to} \tilde{Y}.$$
Looking at Theorem \ref{asip1}(c), we can also conclude that 
$$a_n \max_{k \in K_{n,2}(\cdot)} |T_k|/\sqrt{k} -b_{d,n} \to -\infty \mbox{ a.s.}$$
and the proof of (\ref{DEtheo}) further reduces to showing 
$$a_n \max_{k \in K_n} |T_k|/\sqrt{k} -b_{d,n} \stackrel{d}{\to} \tilde{Y}.$$
Using the same argument as in (\ref{proof1}), we also see that
$$
a_n \max_{1 \le k <k_n} |T_k|/\sqrt{k} -b_{d,n} \to -\infty \mbox{ a.s.}
$$
and we have shown that (\ref{DEtheo}) holds if
$$a_n \max_{1 \le k \le n} |T_k|/\sqrt{k} -b_{d,n} \stackrel{d}{\to} \tilde{Y}.$$
This is the Darling-Erd\H os theorem for normally distributed random vectors which follows from (\ref{Hor}). Thus  (\ref{DEtheo}) has been proven. \\
We now turn to the proof of (\ref{DEtheo1}). By Slutsky's lemma and (\ref{proof1}) it is enough to show that
$$\Delta_n:=a_n \left| \max_{k_n \le k \le n}  |S_k|/\sqrt{k} -  \max_{k_n \le k \le n}  |\Gamma_k^{-1}S_k|/\sqrt{k} \right| \stackrel{\PP}{\to} 0,$$
Using the triangular inequality, it is easy to see that
$$\Delta_n \le a_n \max_{k_n \le k \le n} \left|(I -\Gamma_k)\Gamma_k^{-1}S_k/\sqrt{k}\right|
\le a_n\|I -\Gamma_{k_n}\| \max_{1 \le k \le n} |\Gamma_k^{-1}S_k|/\sqrt{k}.$$
From (\ref{DEtheo}) it follows that $ (\max_{1 \le k \le n} |\Gamma_k^{-1}S_k|/\sqrt{k})/\sqrt{LLn}$ is stochastically bounded. By assumption (\ref{log-log}) we also have that $\|I -\Gamma_{k_n}\| = o((LLn)^{-1})$. Recalling that $a_n=\sqrt{2LLn},$ we see that $\Delta_n \stackrel{\PP}{\to} 0$
and our proof of Theorem \ref{DE-general} is complete. $\Box$
\section{Proof of Theorem \ref{converse}} 
Using the same arguments as in the proof of Theorem \ref{DE-general} we can infer  from (\ref{EXTRAc}) via relations  (\ref{eq-2}) and (\ref{step3f1}) that 
\be \label{M_n}
M_n:=a_n \max_{k_n \le k \le n} |T'_k|/\sqrt{k} - b_{d,n} \stackrel{d}{\to} \tilde{Y}+c,
\ee
where $T'_k =\sum_{j=1}^k \tilde{\Gamma}_j Z_j$ and the random vectors $Z_j$ are i.i.d. with $\mathcal{N}(0, I)$-distribution and $k_n \le \exp((Ln)^{\alpha})$ for some $0 < \alpha <1.$\\[.1cm]
Our first lemma gives an upper bound of $\PP\{M_n > t\}$ via the corresponding probability for the maximum of a subcollection of the random variables $\{|T'_k|/\sqrt{k}: k_n \le k \le n\}$.  (See Lemma 4.3 in \cite{Da-Er} for a related result.)\\[.1cm]
 Let $0 < \xi  < 1$ be fixed. Set 
 $$m_j = [\exp(j\xi/LLn)], j \ge 1 \mbox{ and }N=N_n=[Ln LLn /\xi].$$ Then $m_{N} \le n \le m_{N+1}$. Also note that the sequence $m_j$ depends on $n$ and $\xi$. \\[.1cm]
Next, set 
$$j_n:= \min\{j: m_j \ge Ln\} \mbox{ and }k_n=m_{j_n}$$ so that $j_n \sim \xi^{-1}(LLn)^2$ and $k_n \sim Ln$ as $n \to \infty.$\\ Finally to simplify notation, we set $f_n(y)=(b_{d,n} + y)/a_n, y \in \R$ so that
$$\PP\{M_n > y\}=\PP\left\{\max_{k_n \le k \le n} |T'_k|/\sqrt{k} > f_n(y)\right\}.$$
\begin{lem} \label{lem31}
Given $ 0 <  \delta  < 1$, we have for $y \in \R$ and $n \ge n_0=n_0(\xi,\delta,y),$
$$ (1-\delta)\; \mathbb{P}\left\{ M_n > y+\delta \right\}
\le  \mathbb{P}\left\{ \max_{j_n \le j \le N} | T'_{m_j} | / \sqrt{ m_j }   > f_n(y)  \right\} + 
\PP\{|Z_1| \ge f_n(y)\},$$
provided that $0 < \xi \le \delta^{3}/(36d).$
\end{lem}	
{\bf Proof.} Noting that
\begin{align*}
 &\;\mathbb{P}\left\{ M_n > y+\delta \right\}=\mathbb{P}\left\{ \max_{k_n \le k \le n} |T'_k|/ \sqrt{k}    > f_n(y+ \delta) \right\} \\
 \leq &\;
 	 \mathbb{P}\left\{  \max_{k_n \le k \le n }| T'_k |/ \sqrt{k}    > f_n(y+ \delta), \max_{j_n \le j \le N}  |T'_{m_j} |/ \sqrt{ m_j }   \le f_n(y )  \right\} \\
	 &+  \mathbb{P}\left\{  \max_{j_n \le j \le N} | T'_{m_j}| / \sqrt{ m_j }   >  f_n(y )  \right\} 
,\end{align*}
it is enough to show that 
\begin{eqnarray} \label{eq31}
&&\mathbb{P}\left\{  \max_{k_n \le k \le n} |T'_k |/ \sqrt{k}    > f_n(y+ \delta), \max_{j_n \le j \le N}| T'_{m_j} | /\sqrt{ m_j }  \le  f_n(y )  \right\}\\ \nonumber &\le &\delta\; \mathbb{P}\left\{  \max_{k_n \le k \le n} |T'_k |/ \sqrt{k}    > f_n(y+ \delta)\right\} + \PP\{|Z_1| \ge f_n(y)\} ,
\end{eqnarray}
if $\xi$ is sufficiently small.\\[.2cm]
Consider the following stopping time,
$$\tau := \inf \{ k\ge k_n : | T'_k | / \sqrt{k} > f_n(y+\delta) \}.$$
Then it is obvious that the probability in (\ref{eq31}) is bounded above by
\be \label{eq32}
 \sum_{j=j_n+1}^{N-1 } \sum_{k= m_{j-1}+1}^{ m_j-1}  \mathbb{P}\left\{ \tau = k,  \max_{j_n \le j \le N}| T'_{m_j} | /\sqrt{ m_j }   \le f_n(y )   \right\}  
+ \; \mathbb{P}\left\{ m_{N-1} < \tau \le n \right\} 
\ee
Furthermore, we have for  $j_n +1 \le j \le N-1,$
\begin{align*}
 & \;\sum_{k= m_{j-1}+1}^{ m_j -1}  \mathbb{P}\left\{ \tau = k,  \max_{j_n\le j \le N} | T'_{m_j} | /\sqrt{ m_j }   \le  f_n(y )   \right\}\\
\leq & \;  \sum_{k= m_{j-1}+1}^{ m_j -1} \mathbb{P}\{ \tau = k \} \mathbb{P}\left\{  \frac{ |T'_k - T'_{m_{j+1}}|}  {\sqrt{m_{j+1} - k} }   > \frac{ \sqrt{k} f_n(y+\delta) -\sqrt{m_{j+1} } f_n(y) }{ \sqrt{m_{j+1} - k} }\right\}. 
\end{align*}
Next observe that
\begin{eqnarray*}
&&\max_{m_{j-1}< k \le m_j}  \mathbb{P}\left\{  \frac{ |T'_k - T'_{m_{j+1}} | }{ \sqrt{m_{j+1} - k} }   > \frac{ \sqrt{k} f_n(y+\delta) -\sqrt{m_{j+1} } f_n(y) }{ \sqrt{m_{j+1} - k} }\right\} \\
&\le&  \mathbb{P}\left\{ | Z_1 | > \frac{ \sqrt{m_{j-1} } f_n(y+\delta) -\sqrt{m_{j+1} } f_n(y) }{ \sqrt{m_{j+1} - m_{j-1} } }\right\}. 
\end{eqnarray*}
After some calculation we find that for large enough $n,$
$$  \frac{ \sqrt{m_{j-1} } f_n(y+\delta) -\sqrt{m_{j+1} } f_n(y) }{ \sqrt{m_{j+1} - m_{j-1} } } \ge \frac{\delta}{3\sqrt{\xi}}- 4\sqrt{\xi} \ge \frac{\delta}{6\sqrt{\xi}}$$
where the last inequality holds since $\xi \le  \delta/24.$ 
We trivially have by Markov's inequality, 
$$\PP\{|Z_1| \ge \delta/(6\sqrt{\xi})\} \le 36 \xi \E[|Z_1|^2]/\delta^2= 36\xi d/\delta^2,$$
which is $\le \delta$ by our condition on $\xi.$\\[.2cm]
It follows that
\be \label{eq33}
\sum_{j=j_n+1}^{N-1 } \sum_{k= m_{j-1}+1}^{ m_j -1}  \mathbb{P}\left\{ \tau = k,  \max_{j_n\le j \le N} | T'_{m_j}|/\sqrt{ m_j }   \le f_n(y )   \right\}  \le \delta\; \PP\{k_n \le \tau \le m_{N-1}\}.
\ee
Concerning the second term in (\ref{eq32}) simply note that
\begin{eqnarray*}
&&\mathbb{P}\left\{ m_{N-1} < \tau \le n \right\}\\
&\le & \mathbb{P}\left\{ m_{N-1} < \tau \le n, |T'_{m_{N+1}}|/\sqrt{m_{N+1}} \le f_n(y) \right\} + \PP\{ |T'_{m_{N+1}}|/\sqrt{m_{N+1}} > f_n(y) \}\\
&=&\sum_{k=m_{N-1}+1}^n \mathbb{P}\left\{  \tau =k , |T'_{m_{N+1}}|/\sqrt{m_{N+1}} \le f_n(y) \right\} + \PP\{|Z_1| > f_n(y)\}.
\end{eqnarray*}
Arguing as above, we readily obtain,
\be \label{eq34}
\mathbb{P}\left\{ m_{N-1} < \tau \le n \right\} \le \delta \mathbb{P}\left\{ m_{N-1} < \tau \le n \right\}
+ \PP\{|Z_1| > f_n(y)\}.
\ee
Combining relations (\ref{eq33}) and (\ref{eq34}) and recalling (\ref{eq32}), we see that
\begin{eqnarray*}
&&\mathbb{P}\left\{  \max_{k_n \le k \le n} |T'_k |/\sqrt{k}    > f_n(y+ \delta), \max_{j_n \le j \le N} | T'_{m_j} |/ \sqrt{ m_j }  \le  f_n(y )  \right\} \\
&\le& \delta\;\PP\{k_n \le \tau \le n\} + \PP\{|Z_1| > f_n(y)\}.
\end{eqnarray*}
This implies (\ref{eq31}) since 
$$\PP\{k_n \le \tau \le n\}= \mathbb{P}\left\{ \max_{k_n\le k \le n} |T'_k | / \sqrt{k}  > f_n(y+\delta) \right\},$$
and the proof of Lemma \ref{lem31} is complete. $\Box$\\[.3cm]
We finally need the following lemma,
\begin{lem} \label{lem32}
 Let $Y$ be a $d$-dimensional random vector with distribution $\mathcal{N}(0,\Sigma),$ where $d \ge 2.$
Assume that the largest eigenvalue of $\Sigma$ is equal to 1 and has multiplicity $d-1.$ Denote the remaining (smallest) eigenvalue of $\Sigma$ by $\sigma^2$. Then we have:
$$\PP\{|Y| \ge t\} \le \frac{2}{\sqrt{1-\sigma^2}} \PP\{|Z| \ge t\}, t >0,$$
where $Z: \Omega \to \R^{d-1}$ has a normal$(0, I_{d-1})$-distribution. 
 \end{lem}
{\bf Proof.} If $d \ge 3$, Lemma \ref{lem32} follows  by integrating the inequality given in  Lemma 1(a) of \cite{E-3}.\\[.1cm]
To prove  Lemma \ref{lem32} if $d=2$, we proceed similarly as in \cite{E-3}.
Choose an orthonormal basis $e_1, e_2$ of $\R^2$ consisting of two eigenvectors corresponding to the eigenvalues $1$ and $\sigma^2 \in ]0,1[$ of $\Sigma.$ Then,
$$Y=\sum_{i=1}^2 \langle Y, e_i\rangle e_i =:  \eta_1 e_1 + \sigma \eta_2 e_2 $$
where $\eta_i, 1 \le i \le 2 $ are independent standard normal random variables.\\
It is then obvious that 
$$Y^2 = \eta_1^2 + \sigma^2 \eta_2^2 =: R_1 +R_2,$$
where $R_1$ and $R_2/\sigma^2$ have chi-square distributions with $1$  degree of freedom. Denote the densities of $R_1+R_2$, $R_1$ and $R_2$ by $h, h_1, h_2$.\\
Then  $h_2(y)=h_1(y/\sigma^2)/\sigma^2$ and
$$h(z)=\sigma^{-2}\int_0^z h_1(z-y)h_1(y/\sigma^2)dy, z \ge 0.$$
Using that $h_1(y)=(2\pi)^{-1/2} y^{-1/2} e^{-y/2}, y > 0$, we can infer that
\begin{eqnarray*}
&&h(z)/h_1(z) =\frac{1}{2\pi \sigma} \int_0^z (1-y/z)^{-1/2} y^{-1/2} e^{-(\sigma^{-2} -1)y/2} dy\\
&\le& \frac{1}{\sqrt{2}\pi \sigma} \int_0^{z/2} y^{-1/2} e^{-(\sigma^{-2} -1)y/2} dy 
+  \frac{e^{-(\sigma^{-2} -1)z/4}}{\sqrt{2}\pi \sigma \sqrt{z}}\int_{z/2}^z (1-y/z)^{-1/2}dy\\
&\le&  (2\pi\sigma^2 (\sigma^{-2}-1))^{-1/2} + \sqrt{2z}e^{-(\sigma^{-2} -1)z/4} (\pi \sigma)^{-1}.
\end{eqnarray*}
Employing the trivial inequality $e^{-x/2}\le x^{-1/2}, x > 0,$ it follows that
$$h(z)/h_1(z) \le [(2\pi)^{-1/2} + \sqrt{8}/\pi](1-\sigma^2)^{-1/2} \le 2 (1-\sigma^2)^{-1/2}, z \ge 0.$$
We can conclude that for $t \ge 0,$
$$\PP\{|Y| \ge t\} = \int_{t^2}^{\infty} h(z) dz \le \frac{2}{\sqrt{1-\sigma^2}}\int_{t^2}^{\infty} h_1(z) dz
= \frac{2}{\sqrt{1-\sigma^2}}\PP\{|Z| \ge t\}$$
and Lemma \ref{lem32} has been proven. $\Box$\\[.2cm]
Recall that $d_n=\sqrt{n}/(LLn)^5$ and $\tilde{\Gamma}_n = A(d_n), n \ge 1.$ Let $\{v_1, \ldots, v_d\}$ be an orthonormal basis of $\R^d.$ Then it is easy to see that
\begin{eqnarray*}
\E|X|^2I\{|X|>d_n\} &=&\sum_{i=1}^d \E\langle X, v_i \rangle^2 I\{|X|>d_n\}\\
&=&\sum_{i=1}^d \langle v_i, (I-\tilde{\Gamma}^2_n)v_i\rangle \le d\| I-\tilde{\Gamma}^2_n\|.
\end{eqnarray*}
It is now obvious that  condition (\ref{fconv}) is equivalent to 
$$\|I-\tilde{\Gamma}^2_n \| = O ((LLn)^{-1}) \mbox{ as }n \to \infty,$$
Furthermore, $\|I - \tilde{\Gamma}^2_n \| $ is equal to $1 -\tilde{\lambda}^2_n,$ where $\tilde{\lambda}_n$ is the smallest eigenvalue of $\tilde{\Gamma}_n$ since $\tilde{\Gamma}^2_n$ is symmetric and $I -\tilde{\Gamma}^2_n $ is non-negative definite. So it remains to be shown that (\ref{M_n}) implies
\be \label{eq41}
1 -\tilde{\lambda}^2_n = O((LLn)^{-1}) \mbox{ as }n \to \infty.
\ee
or, equivalently, to show that if (\ref{eq41}) does not hold, we cannot have (\ref{M_n}). \\
To that end we apply Lemma \ref{lem31} with $\delta=1/2$ and we get for $y \in \R,$
\begin{eqnarray}
&&\PP\left\{ M_n \ge y\right\} \nonumber \\
&\le& 2 \mathbb{P}\left\{ \max_{j_n \le j \le N} | T'_{m_j} | /\sqrt{ m_j }   > f_n(y-1/2)  \right\} + 
2 \PP\left\{| Z_1| \ge f_n(y-1/2)\right\} \nonumber\\
&\le& 2 N  \PP\{|\tilde{\Gamma}_n Z_1| \ge f_n(y-1/2\} +  2 \PP\left\{| Z_1| \ge f_n(y-1/2)\right\}. \label{f57}\end{eqnarray}
Here we have used the monotonicity of the sequence $\tilde{\Gamma}_k, k \ge 1$ which implies that $\tilde{\Gamma}_n -\mathrm{Cov}( T'_{m_j}/\sqrt{m_j})$ is non-negative definite for $j_n \le j \le N.$ This allows us to conclude that for $j_n \le j \le N,$
$$\PP\{ |T'_{m_j} | /\sqrt{ m_j }   > x\} \le \PP\{|\tilde{\Gamma}_n Z_1| > x\}, x \in \R.$$
Let $D_n$ be the $d$-dimensional diagonal matrix with  $D_n(i,i)=1, 1 \le i \le d-1$ and $D_n(d,d)=\tilde{\lambda}_n$. Then clearly
$$ \PP\{|\tilde{\Gamma}_{n} Z_1| \ge f_n(y-1/2)\} \le \PP\{|D_{n} Z_1| \ge f_n(y-1/2)\}$$ and we can infer from Lemma \ref{lem32} that
\be \label{eq41a} \PP\{|\tilde{\Gamma}_{n} Z_1| \ge f_n(y-1/2)\} \le 2(1-\tilde{\lambda}_n^2)^{-1/2} \PP\{|Z'| \ge f_n(y-1/2)\},\ee
where $Z'$ is a  $(d-1)$-dimensional normal mean zero random vector with covariance matrix equal to the identity matrix. \\
Using the fact that the square of the Euclidean norm of a $d$-dimensional $\mathcal{N}(0,I)$-distributed random vector $X$ has a gamma distribution with parameters $d/2$ and $2$, one can show that there exist positive constants $C_1(d), C_2(d)$ so that
\be \label{eq42}
C_1(d) t^{d-2}\exp(-t^2/2) \le \PP\{|X| \ge t\} \le C_2(d) t^{d-2}\exp(-t^2/2), t \ge 2d.
\ee
(See Lemma 1 and Lemma 3 in \cite{E-4}, where more precise bounds are given if $d \ge 3$. If $d=1$ this follows directly from   well known bounds  for the tail probabilities of the 1-dimensional normal distribution. If $d=2$ the random variable $|X|^2$ has an exponential distribution and (\ref{eq42}) is trivial.)\\[.1cm]
We can conclude that for large enough $n,$
\begin{eqnarray*}
 \PP\{|Z'| \ge f_n(y-1/2)\} &\le& C_2(d-1)f_n(y-1/2)^{d-3} \exp(-(f_n(y-1/2)^2/2)\\
&\le& C_3(d) f_n(y-1/2)^{-1} \PP\{|Z_1| \ge f_n(y-1/2)\},
\end{eqnarray*}
where we set $C_3(d)= C_2(d-1)/C_1(d).$
Returning to inequality (\ref{eq41a}) and noting that $f_n(y-1/2) \ge \sqrt{\log \log n}$ if $n$ is large, we get in this case,
$$\PP\{|\tilde{\Gamma}_{n} Z_1| \ge f_n(y-1/2)\} \le 2C_3(d)\{(1-\tilde{\lambda}_n^2)\log\log n\}^{-1/2}\PP\{|Z_1| \ge f_n(y-1/2)\}.$$
Applying (\ref{eq42}) once more we find  that 
$$\PP\{|Z_1| \ge f_n(y-1/2)\}=O(N^{-1})=O((\log n \log \log n)^{-1}) \mbox{ as }n \to \infty.$$ 
Recalling (\ref{f57}) we can conclude that if 
$$\limsup_{n \to \infty} (1-\tilde{\lambda}_n^2) \log \log n = \infty,$$ we have for any $y \in \mathbb{R},$
$$\liminf_{n \to \infty} \PP\left\{ M_n > y\right\}=0.$$
Consequently $M_n$  cannot converge in distribution to any variable of the form $\tilde{Y}+c$. $\Box$\\[.2cm]
{\bf Remarks} 
\begin{enumerate}
\item Denote the distribution of $a_n\max_{1 \le k \le n} |S_k|/\sqrt{k} -b_{d,n} $ by $Q_n$. From (\ref{eq-2}) and (\ref{step3f1}) 
 it follows that this sequence is tight if and only if the distributions of $M_n$ form a tight sequence. The above argument actually shows that this last sequence cannot be tight if condition (\ref{fconv}) is not satisfied. Moreover, it is not difficult to prove via Theorem \ref{DE-general} that (\ref{fconv}) implies that the sequence $\{Q_n: n \ge 1\}$  is tight. Thus we have 
$$\{Q_n: n \ge 1\} \mbox{ is tight } \Longleftrightarrow (\ref{fconv}).$$
\item 
Also note that 
$$\PP\left\{a_n \max_{1 \le k \le n} |\Gamma_k^{-1}S_k|/\sqrt{k} - b_n \le 0\right\}
\le \PP\left\{a_n \Lambda_n^{-1}\max_{1 \le k \le n} |S_k|/\sqrt{k} - b_n \le 0\right\},$$
where $\Lambda_n$ is the {\bf largest} eigenvalue of $\Gamma_n$ which in turn is defined as in (\ref{gamma_n}).(Here we can choose any  sequence $c_n$ satisfying condition (\ref{c_n}).)\\ Using this inequality one can show by the same argument as on p. 255 in \cite{E-5} that (\ref{DEtheo1}) implies 
$$1 - \Lambda^2_n = o ((LLn)^{-1}) \mbox{ as }n \to \infty.$$
This is of course weaker than (\ref{log-log}) if $d \ge 2.$
\end{enumerate}\newpage
\section{Some further results}
We first prove the following Darling-Erd\H os type theorem with a shifted limiting distribution.
\begin{theo} \label{shift}
 Let $X, X_n, n \ge 1$ be i.i.d. real-valued  random variables with $\E X^2 =1$ and $\E X=0.$
 Assume that for some $c>0,$
 $$\E X^2I\{|X| \ge t\} \sim c (LLt)^{-1} \mbox{ as }t \to \infty.$$
 Then we have,
 $$a_n \max_{1 \le k \le n} |S_k|/\sqrt{k} - b_{n} \stackrel{d}{\to} \tilde{Y}-c,$$
 where $\tilde{Y}$ and $b_n$ are defined  as in (\ref{1956}).
\end{theo}
{\bf Proof.} (i) Set $\sigma_n^2 =\E X^2I\{|X| \le \sqrt{n}\}$ and let $1 \le k_n \le \exp((Ln)^{\alpha})$ for some $0 < \alpha <1.$  Then we have by Theorem \ref{DE-general} and the argument in (\ref{proof1}),
$$a_n \max_{k_n \le k \le n} |S_k|/ \sqrt{k}\sigma_k  - b_{n} \stackrel{d}{\to} \tilde{Y},$$
which trivially implies for any sequence $ \rho_n$ of positive real numbers converging to 1,
\be \label{f61}
\rho_n a_n \max_{k_n \le k \le n} |S_k|/\sqrt{k}\sigma_k - \rho_n b_{n} \stackrel{d}{\to} \tilde{Y}
\ee
Set $k_n = [\exp((Ln)^{\alpha}]$, where $0< \alpha <1.$ Then it is easy to see that
\begin{eqnarray*}
&&\PP\left\{a_n \max_{1 \le k \le n} \frac{|S_k|}{\sqrt{k}} - b_n \le y\right\}\\
&\le& \PP\left\{\sigma_{k_n} a_n \max_{k_n \le k \le n} \frac{|S_k|}{\sqrt{k}\sigma_k} - b_n \le y\right\}\\&=&\PP\left\{\sigma_{k_n} a_n \max_{k_n \le k \le n} \frac{|S_k|}{\sqrt{k}\sigma_k} - \sigma_{k_n} b_n \le y +(1-\sigma_{k_n})b_n\right\}
\end{eqnarray*}
Noticing that $(1-\sigma_{k_n})b_n \sim (1 -\sigma^2_{k_n})(2LLn)/(1+\sigma_{k_n})\sim c(LLk_n)^{-1}LLn$ (since $\sigma^2_{k_n} \to \E X^2=1$), it is clear that $(1-\sigma_{k_n})b_n \to c/\alpha$ as $n \to \infty.$\\[.1cm]
By (\ref{f61}) (with $\rho_n=\sigma_{k_n}$) this last sequence of probabilities converges to\\ $\PP\{\tilde{Y} \le y +c/\alpha\}.$\\
Since this holds for any $0 < \alpha <1,$ it follows that
\be \label{f62}
\limsup_{n \to \infty} \PP\left\{a_n \max_{1 \le k \le n} \frac{|S_k|}{\sqrt{k}} - b_n \le y\right\} \le 
\PP\{\tilde{Y} \le y +c\}.
\ee
(ii) Similarly, we have,
\begin{eqnarray*}
&&\PP\left\{a_n \max_{1 \le k \le n} \frac{|S_k|}{\sqrt{k}} - b_n \le y\right\}\\
&\ge& \PP\left\{\sigma_{n} a_n \max_{1 \le k \le n} \frac{|S_k|}{\sqrt{k}\sigma_k} - \sigma_{n} b_n \le y +(1-\sigma_{n})b_n\right\},
\end{eqnarray*}
where $(1-\sigma_n)b_n \to c$ as $n \to \infty.$\\
Applying (\ref{f61}) (with $k_n=1$ and $\rho_n=\sigma_n$), we obtain that
$$\liminf_{n \to \infty} \PP\left\{a_n \max_{1 \le k \le n} \frac{|S_k|}{\sqrt{k}} - b_n \le y\right\} \ge 
\PP\{\tilde{Y} \le y +c\}$$
and Theorem \ref{shift} has been proven. $\Box$\\[.2cm]
We finally mention the following result for real-valued random variables given in \cite{kls} where it is shown that if $\E X=0, \E X^2 =1$ and $\E X^2 LL|X| < \infty$, then one has
\be \label{question}
2LLn \left(\sup_{k \ge n} \frac{|S_k|}{\sqrt{2kLLk}}-1\right) -\frac{3}{2} LLLn + LLLLn + \log(3/\sqrt{8}) \stackrel{d}{\to} \tilde{Y}.
\ee
The authors asked whether this result can hold under the finite second moment assumption.\\
 Using Theorem \ref{asip1} in combination with Theorem 1.1 in \cite{kls}, we  obtain the following general result:
$$2LLn \left(\sup_{k \ge n} \frac{|S_k|}{\sqrt{2 k LLk}\sigma_k}-1\right) -\frac{3}{2} LLLn + LLLLn + \log(3/\sqrt{8}) \stackrel{d}{\to} \tilde{Y},$$
where $\sigma_n^2 = \E X^2 I\{|X| \le c_n\}$ and $c_n$ is a non-decreasing sequence of positive real numbers satisfying condition (\ref{c_n}). As in \cite{E-5} this implies that (\ref{question}) holds if and only if condition (\ref{log-log}) is satisfied.
{\bf Acknowledgements}\\
The authors would like to thank the referee for some useful suggestions helping us to improve the presentation of our results.


\begin{thebibliography}{xx}
\bibitem{bai} Bai, Z. D. (1989). A theorem of Feller revisited. {\it Ann. Probab.} {\bf 17}, 385--395.
\bibitem{be}Bertoin, J. (1998). Darling-Erd\H os theorems for normalized sums of i.i.d. variables close to a stable law. {\it Ann. Probab.} {\bf 26}, 832--852.
\bibitem{Bha} Bhatia, R. (1997). {\it Matrix Analysis}. Springer. New York.
\bibitem{csw} Cs\"org\H o, M.; Szyskowicz, B. and Wang, Q. (2003) Darling-Erd\H os theorem for self-normalized sums. {\it Ann. Probab.} {\bf 31}, 676--692.
\bibitem {Da-Er} Darling, D. A. and Erd\H os, P. (1956). A limit theorem for the maximum of normalized sums of independent random variables. {\it Duke Math. J.} {\bf 23} 143--155.
\bibitem{E-5} Einmahl, U. (1989). The Darling-Erd\H{o}s theorem for sums of i.i.d. random variables. {\it Probab. Th. Rel. Fields} \textbf{82}, 241--257.
\bibitem{E-3} Einmahl, U. (1991). On the almost sure behavior of sums of iid random variables in Hilbert space. {\it Ann. Probab.} {\bf 19}, 1227--1263.
\bibitem{E-4} Einmahl, U. (1992). Exact convergence rates for the bounded law of the iterated logarithm in Hilbert space. {\it Probab. Th. Rel. Fields} {\bf 92}, 177--194.
\bibitem{E-2} Einmahl, U. (2009).  A new strong invariance principle for sums of independent random vectors. {\it Journal of Mathematical Sciences} {\bf163}, 311--327. 
\bibitem{EM} Einmahl, U. and Mason, D. M. (1989). Darling-Erd\H os theorems for martingales. {\it J. Theoret. Probab.} {\bf 2}, 437--460.
\bibitem{Fe} Feller, W. (1946). The law of the iterated logarithm for identically distributed random variables. {\it Ann. Math.} {\bf 47}, 631--638.
\bibitem{Hor} Horv\'ath, L. (1994). Likelihood method for testing changes in the parameters of normal observations. {\it Ann. Stat.}   \textbf{21}(2), 671--680.
\bibitem{kls} Khoshnevisan, D.; Levin, D. A. and Shi, Z. (2005). An extreme-value analysis of the LIL for Brownian motion. {\it Elect. Comm. in Probab.} {\bf 10}, 196--206.
\bibitem{ood} Oodaira, H. (1976). Some limit theorems for the maximum of normalized sums of weakly dependent random variables. In:  Proceedings of the third Japan-USSR symposium on Probability Theory.
{\it Lecture Notes in Mathematics} {\bf 550}, 1--13. Springer, Berlin.
\bibitem{sakh} Sakhanenko, A.I. (2000). 
A new way to obtain estimates in the invariance principle. In:  Proceedings of High Dimensional Probability II, Birkh\"auser, {\it Progress in Probability}, \textbf{47}, 223--245.
\bibitem{sho} Shorack, G. (1979). Extension of the Darling and Erd\H os theorem on the maximum of normalized sums. {\it Ann. Probab.} {\bf 7}, 1092--1096.
\end{thebibliography}


\end{document}